\documentstyle[11pt]{article}
\begin{document}

%
%
%
\newtheorem{theorem}      {Th\'eor\`eme}[section]
\newtheorem{theorem*}     {theorem}
\newtheorem{proposition}  [theorem]{Proposition}
\newtheorem{definition}   [theorem]{Definition}
\newtheorem{e-lemme}        [theorem]{Lemma}
\newtheorem{cor}   [theorem]{Corollaire}
\newtheorem{resultat}     [theorem]{R\'esultat}
\newtheorem{eexercice}    [theorem]{Exercice}
\newtheorem{rrem}    [theorem]{Remarque}
\newtheorem{pprobleme}    [theorem]{Probl\`eme}
\newtheorem{eexemple}     [theorem]{Exemple}
\newcommand{\preuve}      {\paragraph{Preuve}}
\newenvironment{probleme} {\begin{pprobleme}\rm}{\end{pprobleme}}
\newenvironment{remarque} {\begin{rremarque}\rm}{\end{rremarque}}
\newenvironment{exercice} {\begin{eexercice}\rm}{\end{eexercice}}
\newenvironment{exemple}  {\begin{eexemple}\rm}{\end{eexemple}}
%
%
\newtheorem{e-theo}      [theorem]{Theorem}
\newtheorem{theo*}     [theorem]{Theorem}
\newtheorem{e-pro}  [theorem]{Proposition}
\newtheorem{e-def}   [theorem]{Definition}
\newtheorem{e-lem}        [theorem]{Lemma}
\newtheorem{e-cor}   [theorem]{Corollary}
\newtheorem{e-resultat}     [theorem]{Result}
\newtheorem{ex}    [theorem]{Exercise}
\newtheorem{e-rem}    [theorem]{Remark}
\newtheorem{prob}    [theorem]{Problem}
\newtheorem{example}     [theorem]{Example}
\newcommand{\proof}         {\paragraph{Proof~: }}
\newcommand{\hint}          {\paragraph{Hint}}
\newcommand{\heuristicproof}{\paragraph{heuristic proof}}
\newenvironment{e-probleme} {\begin{e-pprobleme}\rm}{\end{e-pprobleme}}
\newenvironment{e-remarque} {\begin{e-rremarque}\rm}{\end{e-rremarque}}
\newenvironment{e-exercice} {\begin{e-eexercice}\rm}{\end{e-eexercice}}
\newenvironment{e-exemple}  {\begin{e-eexemple}\rm}{\end{e-eexemple}}
\newcommand{\reell}    {{{\rm I\! R}^l}}
\newcommand{\reeln}    {{{\rm I\! R}^n}}
\newcommand{\reelk}    {{{\rm I\! R}^k}}
\newcommand{\reelm}    {{{\rm I\! R}^m}}
\newcommand{\reelp}    {{{\rm I\! R}^p}}
\newcommand{\reeld}    {{{\rm I\! R}^d}}
\newcommand{\reeldd}   {{{\rm I\! R}^{d\times d}}}
\newcommand{\reelnn}   {{{\rm I\! R}^{n\times n}}}
\newcommand{\reelnd}   {{{\rm I\! R}^{n\times d}}}
\newcommand{\reeldn}   {{{\rm I\! R}^{d\times n}}}
\newcommand{\reelkd}   {{{\rm I\! R}^{k\times d}}}
\newcommand{\reelkl}   {{{\rm I\! R}^{k\times l}}}
\newcommand{\reelN}    {{{\rm I\! R}^N}}
\newcommand{\reelM}    {{{\rm I\! R}^M}}
\newcommand{\reelplus} {{{\rm I\! R}^+}}
\newcommand{\reelo}    {{{\rm I\! R}\setminus\{0\}}}
\newcommand{\reld}    {{{\rm I\! R}_d}}
\newcommand{\relplus} {{{\rm I\! R}_+}}
\newcommand{\1}        {{\bf 1}}

\newcommand{\cov}      {{\hbox{cov}}}
\newcommand{\sss}      {{\cal S}}
\newcommand{\indic}    {{{\rm I\!\! I}}}
\newcommand{\pp}       {{{\rm I\!\!\! P}}}
\newcommand{\qq}       {{{\rm I\!\!\! Q}}}
\newcommand{\ee}       {{{\rm I\! E}}}

\newcommand{\B}        {{{\rm I\! B}}}
\newcommand{\cc}       {{{\rm I\!\!\! C}}}
\newcommand{\HHH}        {{{\rm I\! H}}}
\newcommand{\N}        {{{\rm I\! N}}}
\newcommand{\R}        {{{\rm I\! R}}}
\newcommand{\D}        {{{\rm I\! D}}}
\newcommand{\Z}       {{{\rm Z\!\! Z}}}
\newcommand{\C}        {{\bf C}}	
\newcommand{\T}        {{\bf T}}	
\newcommand{\E}        {{\bf E}}	
\newcommand{\rfr}[1]    {\stackrel{\circ}{#1}}
\newcommand{\equiva}    {\displaystyle\mathop{\simeq}}
\newcommand{\eqdef}     {\stackrel{\triangle}{=}}
\newcommand{\limps}     {\mathop{\hbox{\rm lim--p.s.}}}
\newcommand{\Limsup}    {\mathop{\overline{\rm lim}}}
\newcommand{\Liminf}    {\mathop{\underline{\rm lim}}}
\newcommand{\Inf}       {\mathop{\rm Inf}}
\newcommand{\vers}      {\mathop{\;{\rightarrow}\;}}
\newcommand{\versup}    {\mathop{\;{\nearrow}\;}}
\newcommand{\versdown}  {\mathop{\;{\searrow}\;}}
\newcommand{\vvers}     {\mathop{\;{\longrightarrow}\;}}
\newcommand{\cvetroite} {\mathop{\;{\Longrightarrow}\;}}
\newcommand{\ieme}      {\hbox{i}^{\hbox{\smalltype\`eme}}}
\newcommand{\eqps}      {\, \buildrel \rm \hbox{\rm\smalltype p.s.} \over = \,}
\newcommand{\eqas}      {\,\buildrel\rm\hbox{\rm\smalltype a.s.} \over = \,}
\newcommand{\argmax}    {\hbox{{\rm Arg}}\max}
\newcommand{\argmin}    {\hbox{{\rm Arg}}\min}
\newcommand{\indep}{\perp\!\!\!\!\perp}
\newcommand{\abs}[1]{\left| #1 \right|}
\newcommand{\crochet}[2]{\langle #1 \,,\, #2 \rangle}
\newcommand{\espc}[3]   {E_{#1}\left(\left. #2 \right| #3 \right)}
\newcommand{\rang}{\hbox{rang}}
\newcommand{\rank}{\hbox{rank}}
\newcommand{\signe}{\hbox{signe}}
\newcommand{\sign}{\hbox{sign}}

\newcommand\hA{{\widehat A}}
\newcommand\hB{{\widehat B}}
\newcommand\hC{{\widehat C}}
\newcommand\hD{{\widehat D}}
\newcommand\hE{{\widehat E}}
\newcommand\hF{{\widehat F}}
\newcommand\hG{{\widehat G}}
\newcommand\hH{{\widehat H}}
\newcommand\hI{{\widehat I}}
\newcommand\hJ{{\widehat J}}
\newcommand\hK{{\widehat K}}
\newcommand\hL{{\widehat L}}
\newcommand\hM{{\widehat M}}
\newcommand\hN{{\widehat N}}
\newcommand\hO{{\widehat O}}
\newcommand\hP{{\widehat P}}
\newcommand\hQ{{\widehat Q}}
\newcommand\hR{{\widehat R}}
\newcommand\hS{{\widehat S}}
\newcommand\hTT{{\widehat T}}
\newcommand\hU{{\widehat U}}
\newcommand\hV{{\widehat V}}
\newcommand\hW{{\widehat W}}
\newcommand\hX{{\widehat X}}
\newcommand\hY{{\widehat Y}}
\newcommand\hZ{{\widehat Z}}

\newcommand\ha{{\widehat a}}
\newcommand\hb{{\widehat b}}
\newcommand\hc{{\widehat c}}
\newcommand\hd{{\widehat d}}
\newcommand\he{{\widehat e}}
\newcommand\hf{{\widehat f}}
\newcommand\hg{{\widehat g}}
\newcommand\hh{{\widehat h}}
\newcommand\hi{{\widehat i}}
\newcommand\hj{{\widehat j}}
\newcommand\hk{{\widehat k}}
\newcommand\hl{{\widehat l}}
\newcommand\hm{{\widehat m}}
\newcommand\hn{{\widehat n}}
\newcommand\ho{{\widehat o}}
\newcommand\hp{{\widehat p}}
\newcommand\hq{{\widehat q}}
\newcommand\hr{{\widehat r}}
\newcommand\hs{{\widehat s}}
\newcommand\htt{{\widehat t}}
\newcommand\hu{{\widehat u}}
\newcommand\hv{{\widehat v}}
\newcommand\hw{{\widehat w}}
\newcommand\hx{{\widehat x}}
\newcommand\hy{{\widehat y}}
\newcommand\hz{{\widehat z}}

\newcommand\tA{{\widetilde A}}
\newcommand\tB{{\widetilde B}}
\newcommand\tC{{\widetilde C}}
\newcommand\tD{{\widetilde D}}
\newcommand\tE{{\widetilde E}}
\newcommand\tF{{\widetilde F}}
\newcommand\tG{{\widetilde G}}
\newcommand\tH{{\widetilde H}}
\newcommand\tI{{\widetilde I}}
\newcommand\tJ{{\widetilde J}}
\newcommand\tK{{\widetilde K}}
\newcommand\tL{{\widetilde L}}
\newcommand\tM{{\widetilde M}}
\newcommand\tN{{\widetilde N}}
\newcommand\tOO{{\widetilde O}}
\newcommand\tP{{\widetilde P}}
\newcommand\tQ{{\widetilde Q}}
\newcommand\tR{{\widetilde R}}
\newcommand\tS{{\widetilde S}}
\newcommand\tTT{{\widetilde T}}
\newcommand\tU{{\widetilde U}}
\newcommand\tV{{\widetilde V}}
\newcommand\tW{{\widetilde W}}
\newcommand\tX{{\widetilde X}}
\newcommand\tY{{\widetilde Y}}
\newcommand\tZ{{\widetilde Z}}

\newcommand\ta{{\widetilde a}}
\newcommand\tb{{\widetilde b}}
\newcommand\tc{{\widetilde c}}
\newcommand\td{{\widetilde d}}
\newcommand\te{{\widetilde e}}
\newcommand\tf{{\widetilde f}}
\newcommand\tg{{\widetilde g}}
\newcommand\th{{\widetilde h}}
\newcommand\ti{{\widetilde i}}
\newcommand\tj{{\widetilde j}}
\newcommand\tk{{\widetilde k}}
\newcommand\tl{{\widetilde l}}
\newcommand\tm{{\widetilde m}}
\newcommand\tn{{\widetilde n}}
\newcommand\tio{{\widetilde o}}
\newcommand\tp{{\widetilde p}}
\newcommand\tq{{\widetilde q}}
\newcommand\tr{{\widetilde r}}
\newcommand\ts{{\widetilde s}}
\newcommand\tit{{\widetilde t}}
\newcommand\tu{{\widetilde u}}
\newcommand\tv{{\widetilde v}}
\newcommand\tw{{\widetilde w}}
\newcommand\tx{{\widetilde x}}
\newcommand\ty{{\widetilde y}}
\newcommand\tz{{\widetilde z}}

\newcommand\bA{{\overline A}}
\newcommand\bB{{\overline B}}
\newcommand\bC{{\overline C}}
\newcommand\bD{{\overline D}}
\newcommand\bE{{\overline E}}
\newcommand\bFF{{\overline F}}
\newcommand\bG{{\overline G}}
\newcommand\bH{{\overline H}}
\newcommand\bI{{\overline I}}
\newcommand\bJ{{\overline J}}
\newcommand\bK{{\overline K}}
\newcommand\bL{{\overline L}}
\newcommand\bM{{\overline M}}
\newcommand\bN{{\overline N}}
\newcommand\bO{{\overline O}}
\newcommand\bP{{\overline P}}
\newcommand\bQ{{\overline Q}}
\newcommand\bR{{\overline R}}
\newcommand\bS{{\overline S}}
\newcommand\bT{{\overline T}}
\newcommand\bU{{\overline U}}
\newcommand\bV{{\overline V}}
\newcommand\bW{{\overline W}}
\newcommand\bX{{\overline X}}
\newcommand\bY{{\overline Y}}
\newcommand\bZ{{\overline Z}}

\newcommand\ba{{\overline a}}
\newcommand\bb{{\overline b}}
\newcommand\bc{{\overline c}}
\newcommand\bd{{\overline d}}
\newcommand\be{{\overline e}}
\newcommand\bff{{\overline f}}
\newcommand\bg{{\overline g}}
\newcommand\bh{{\overline h}}
\newcommand\bi{{\overline i}}
\newcommand\bj{{\overline j}}
\newcommand\bk{{\overline k}}
\newcommand\bl{{\overline l}}
\newcommand\bm{{\overline m}}
\newcommand\bn{{\overline n}}
\newcommand\bo{{\overline o}}
\newcommand\bp{{\overline p}}
\newcommand\bq{{\overline q}}
\newcommand\br{{\overline r}}
\newcommand\bs{{\overline s}}
\newcommand\bt{{\overline t}}
\newcommand\bu{{\overline u}}
\newcommand\bv{{\overline v}}
\newcommand\bw{{\overline w}}
\newcommand\bx{{\overline x}}
\newcommand\by{{\overline y}}
\newcommand\bz{{\overline z}}

%
\newcommand{\AAA}{{\cal A}}
\newcommand{\BB}{{\cal B}}
\newcommand{\CC}{{\cal C}}
\newcommand{\DD}{{\cal D}}
\newcommand{\EE}{{\cal E}}
\newcommand{\FF}{{\cal F}}
\newcommand{\GG}{{\cal G}}
\newcommand{\HH}{{\cal H}}
\newcommand{\II}{{\cal I}}
\newcommand{\JJ}{{\cal J}}
\newcommand{\KK}{{\cal K}}
\newcommand{\LL}{{\cal L}}
\newcommand{\NN}{{\cal N}}
\newcommand{\MM}{{\cal M}}
\newcommand{\OO}{{\cal O}}
\newcommand{\PP}{{\cal P}}
\newcommand{\QQ}{{\cal Q}}
\newcommand{\RR}{{\cal R}}
\newcommand{\SS}{{\cal S}}
\newcommand{\TT}{{\cal T}}
\newcommand{\UU}{{\cal U}}
\newcommand{\VV}{{\cal V}}
\newcommand{\WW}{{\cal W}}
\newcommand{\XX}{{\cal X}}
\newcommand{\YY}{{\cal Y}}
\newcommand{\ZZ}{{\cal Z}}
\newcommand{\tbullet}{$\bullet$}
\newcommand{\ot}{\leftarrow}
\newcommand{\newblock}{}
\newcommand{\carre}{\hfill$\Box$}
\newcommand{\carreb}{\hfill\rule{0.25cm}{0.25cm}}
%
%
\newcommand{\dontforget}[1]
{{\mbox{}\\\noindent\rule{1cm}{2mm}\hfill don't forget : #1 \hfill\rule{1cm}{2mm}}\typeout{---------- don't forget : #1 ------------}}
\newcommand{\note}[2]
{ \noindent{\sf #1 \hfill \today}

\noindent\mbox{}\hrulefill\mbox{}
\begin{quote}\begin{quote}\sf #2\end{quote}\end{quote}
\noindent\mbox{}\hrulefill\mbox{}
\vspace{1cm}
}
\newcommand{\rond}[1]     {\stackrel{\circ}{#1}}
\newcommand{\rondf}       {\stackrel{\circ}{\FF}}
\newcommand{\point}[1]     {\stackrel{\cdot}{#1}}

\newcommand\relatif{{\rm \rlap Z\kern 3pt Z}}

\def\diagram#1{\def\normalbaselines{\baselineskip=Opt
\lineskip=10Opt\lineskiplimit=1pt}  \matrix{#1}}

{\huge Segre varieties  and Lie symmetries    }

\bigskip 

{\bf  Alexandre SUKHOV}

\bigskip

{\footnotesize
Univesrit\'e des Sciences et Technologies de Lille, Laboratoire d'Arithm\'etique
 - G\'eom\'etrie - Analyse - Topologie, Unit\'e Mixte de Recherche 8524, U.F.R. de 
Math\'ematique, 59655 Villeneuve d'Ascq, Cedex, France.

e-mail: Alexandre.Sukhov@agat.univ-lille1.fr}

\bigskip

\bigskip

{\bf Abstract.}  We show that biholomorphic automorphisms of a real analytic hypersurface in $\cc^{n+1}$ 
can be considered as (pointwise) Lie  symmetries of a  
holomorphic completely overdetermined involutive second order PDE system 
defining its Segre family.  Using the classical 
S.Lie method we obtain a complete description of infinitesimal symmetries of such a 
system and give a new proof of some well known results of CR geometry.

\section{Introduction and result}

In the present paper we apply the Lie method of studying of 
PDE symmetries to a special but geometrically important class of holomorphic 
completely overdetermined 
second order PDE systems, i.e. systems of the form

\begin{eqnarray}
\label{1}
({\cal S}): u^k_{x_ix_j} = F^k_{ij}(x,u,u_x), k = 1,...,m, i,j = 1,...,n
\end{eqnarray}
where $x = (x_1,...,x_n)$ are independent variables, $u(x) = (u^1(x),...,u^m(x))$ 
are unknown holomorphic functions (dependent variables), $u_x = (u^j_{x_i})$ and 
$F^k_{ij}= F^k_{ji}$ are holomorphic functions.

Denote by $J^{r}_{n,m}$ the r-jet space \cite{BlKu}, \cite{Ol} of r-jets of holomorphic 
maps from $\cc^n$ to $\cc^m$. Set $u^{(1)} = (u^1_1,...,u^1_n,...,u^{m}_1,...,u^{m}_n)$,...
, $u^{(s)} = (u^j_{\tau})$ with $j= 1,...,m$, $\tau = (\tau_1,...\tau_s)$, 
$\tau_1\leq \tau_2 \leq ...\leq \tau_s$ and use it as the {\it natural coordinates} 
$(x,u,u^{(1)},...,u^{(r)})$ on
the jet space:  for every holomorphic (near a point $p$) function 
$u = f(x): \cc^n \longrightarrow \cc^m$ the natural coordinates of the 
corresponding $r$-jet $j^r_p(f) \in J^r_{n,m}$ defining by $f$ at $p$ are  
$x_j = p_j$, $u^k =f^k(p)$, $u^{k}_{\tau_1...,\tau_s} = 
\frac{\partial^{s} f^k(p)}{\partial x_{\tau_1}...\partial x_{\tau_s}}$.

We say that  a system (\ref{1}) is {\it involutive} if the  differential forms
 $du^k_i - \sum_j F^{k}_{ij}(x,u,u^{(1)})dx_j$, $du^k - \sum_i u^k_idx_i$ define a 
completely integrable distribution on the tangent 
bundle $T(J^1_{n,m})$,   i.e. satisfy  the Frobenius involutivity condition.

 Solutions of such a system are holomorphic vector valued functions $u = u(x)$; 
denote by $\Gamma_u$ the graph of a solution $u$. A symmetry group $Sym(S)$ 
(see for instance \cite{BlKu}, \cite{Ol}) of a system  is a local complex transformation 
group $G$ acting on a domain in the space $\cc^n_x \times \cc^m_u$ of independent and dependent variables 
with the following property: for every solution $u(x)$ of $({\cal S})$ and every $g \in G$ such
that the image $g(\Gamma_u)$ is defined, it is a graph of a solution of $({\cal S})$.
Sometimes the {\it largest} symmetry group $({\cal S})$ is of main interest (and so 
we write {\it the} symmetry group); for us this is not very essential since our methods give a description of ${\it any}$ symmetry group for 
given system.  A holomorphic vector field 

\begin{eqnarray}
\label{2}
X = \sum_j \theta_j(x,u) \frac{\partial}{\partial x_j} + \sum_{\mu}\eta^{\mu}
\frac{\partial}{\partial u^{\mu}} 
\end{eqnarray}
generating  a complex one-parameter group of symmetries of a system of PDE $({\cal S})$ is called 
{\it an infinitesimal symmetry} of this system. All these fields form a complex Lie algebra 
with respect to the Lie bracket which  is  denoted by $Lie({\cal S})$.

 Using the notation $w = (x,u)$, we fix a point $(x,u)$ and set 

\begin{eqnarray*}
& &\alpha_j(x,u) = (\theta_{j_{w_1}}(x,u),...,\theta_{j_{w_{n+m}}}(x,u)), 
\alpha(x,u) = (\alpha_1,...,\alpha_n),\\
& &\beta^k(x,u) = 
(\eta^k_{w_1}(x,u),...,\eta^k_{w_{n+m}}(x,u)), \beta(x,u) = (\beta^1,...,\beta^m),\\
& &\gamma(x,u) = (\theta_{1_{x_1w_1}}(x,u),..., \theta_{1_{x_1w_{n+m}}}(x,u)),\\
& &\delta(x,u) = (\eta^1(x,u),...,\eta^m(x,u)),\\ 
& &\varepsilon(x,u) = 
(\theta_1(x,u),...,\theta_n(x,u))
\end{eqnarray*}

 We call the vector   
\begin{eqnarray*}
\omega(x,u) = (\alpha(x,u),\beta(x,u),\gamma(x,u),\delta(x,u),\varepsilon(x,u))
\end{eqnarray*}
of $\cc^{(n+m+2)(n+m)}$ the {\it initial date} of an infinitesimal symmetry $X$ 
of the form (\ref{2}) at the point $(x,u)$.

Our main result is the following

\begin{e-theo}
\label{theo3.2}
Let $({\cal S})$ be a holomorphic completely overdetermined second order involutive 
system with $n$ independent and $m$ dependent variables. Then the Taylor expansions 
of coefficients of any infinitesimal symmetry 
$X \in Lie({\cal S})$ at a fixed point $(x,u)$ are uniquely determined by the 
initial date $\omega(x,u)$ that is the linear map $Lie({\cal S}) \longrightarrow 
\cc^{(n+m+2)(n+m)}$ defined by $X \mapsto \omega(x,u)$ is injective.
\end{e-theo}

 In the special case $n = 1$, $ m = 1$ i.e. in the case of ordinary second order differential 
equation this result was obtained by A.Tresse \cite{Tr} (a student of S.Lie) in 1896: he 
proved that the symmetry group of a second order differential equation is a Lie group 
of dimension $\leq 8$ (as it was observed by B.Segre, this implies that the automorphism group 
of a real analytic Levi nondegenerate hypersurface in $\cc^2$ is a finite dimensional 
real Lie group). A very clear proof of the same fact is contained in the 
known paper of L.E. Dickson \cite{Di} (another former student of S.Lie) 
inspired by  the lectures of S.Lie at the end of the XIX century. 
So it is quite possible that the basic idea goes back to S.Lie himself. 
Our proof is based on the direct generalization  of the Lie - Tresse - Dickson method.

We point out that this method is quite elementary and constructive and gives 
an efficient  recursive algorithm 
for  determination of the Taylor expansions of the coefficients of an infinitesimal
symmetry at a given point. The symmetry group of $({\cal S})$ then can be parametrized 
 by the exponential map (in a suitable neighborhood of the identity). If a point $(x,u)$ is fixed, the components of the 
initial date $\omega(x,u)$ are local parameters for the symmetry group. 
The  number of these parameters is equal to $(n + m + 2)(n+ m)$; however, {\it in general 
they  are not independent}. In our proof of theorem we consider only those equations on 
the coefficients of an infinitesimal symmetry which are necessary in order to conclude. In general, additional equations can occur. So in the general case the parameters may satisfy 
some additional relations and the actual dimension of a symmetry group may be smaller
than $(n + m + 2)(n +m)$. 

\begin{e-cor} 
A symmetry group of a holomorphic completely overdetermined involutive 
second order system with $n$ independent and $m$ dependent variables is a 
local complex Lie transformation group of dimension $\leq (n + m + 2)(n + m)$.
\end{e-cor}
 
We will show that this estimate is precise. In the special case $m = 1$ this last result can also be deduced from Chern's 
solution of the equivalence problem for completely overdetermined involutive 
second order systems with one dependent variable \cite{Ch}.

In the next section we consider applications of  theorem \ref{theo3.2} to CR geometry.

\section{Segre varieties, holomorphic maps and PDE symmetries}

Let $\Gamma$ be a real analytic {\it Levi nondegenerate} hypersurface in $\cc^{n+1}$
 and let $Aut(\Gamma)$ denote its biholomorphism group (all our considerations are 
local). Denote by $Z = (z,w) \in \cc^n \times \cc$ the standard coordinates in $\cc^{n+1}$. 
For a fixed point $\zeta \in \cc^{n+1}$ close enough to $\Gamma$ consider the 
{\it complex} hypersurface $Q(\zeta) = \{ Z: r(Z,\zeta) = 0 \}$. It is called 
the {\it Segre variety} ( \cite{Se}). 
The basic property of the Segre varieties (\cite{DW},\cite{DF},\cite{We}) 
is their biholomorphic invariance: for every automorphism 
$f \in Aut(\Gamma)$ and any $\zeta$ one has $f(Q(\zeta)) = Q(\overline{f}(\overline\zeta))$. 
   B.Segre observed that for $n= 1$ i.e. in $\cc^2$ the set of Segre varieties of 
$\Gamma$ (which is called {\it the Segre family} of $\Gamma$) is a 
regular two parameter family of holomorphic curves and so represents the trajectories 
of solutions of a holomorphic second order ordinary differential equation (see also 
\cite{Ca}, \cite{Ch},\cite{We}). This important observation can be generalized as follows.

After a biholomorphic change of coordinates in a neighborhood of the 
origin $\Gamma$ is given by the equation $\{ w + \overline{w} + \sum_{j=1}^n \varepsilon_j z_j\overline{z}_j 
+ R(Z,\overline{Z}) = 0 \}$ where $\varepsilon_j = 1$ or $-1$ and $R = o(\vert Z \vert^2)$.
For every point $\zeta$ the corresponding Segre variety 
is given by $w + \zeta_{n+1} + \sum_{j=1}^n \varepsilon_jz_j\zeta_j + R(Z,\zeta) = 0$. 
If we  consider the variables $x_j =z_j$ as independent ones 
and the variable $w = u(x)$ as dependent one, then this equation can be rewritten in the 
form

\begin{eqnarray}
\label{3}
 u + \zeta_{n+1} + \sum_{j=1}^n \varepsilon_jx_j\zeta_j + R((x,u),\zeta) = 0
\end{eqnarray}

 Taking the derivatives in $x_k$ we obtain the equations 

\begin{eqnarray}
\label{4}
u_{x_k} + \varepsilon_k\zeta_k + R_{x_k}(x,u,\zeta) + R_u(x,u,\zeta)u_{x_k} = 0, 
k = 1,...,n 
\end{eqnarray}

The equations (\ref{3}), (\ref{4}) and the implicit function theorem 
imply that  $\zeta$    is an analytic function 
$\zeta= \zeta(x,u,u_{x_1},...,u_{x_n})$; taking again the partial derivatives in 
$x_j$ in (\ref{4})  
and using the obtained  expression for $\zeta$ in order to eliminate it from the equations,
 we obtain that every $u$ given by (\ref{3}) is a solution of a  completely overdermined second order holomorphic PDE  system 
$({\cal S}_{\Gamma})$ of the form (\ref{1}).   This  system is necessarily 
involutive since     the family of its solutions (\ref{3}) define a completely 
integrable distribution on the tangent space $T(J^1_{n,1})$. 
  The property of biholomorphic invariance of the 
Segre varieties means that any biholomorphism of $\Gamma$ 
transforms the graph of a solution of $({\cal S}_{\Gamma})$ to the graph of another solution, 
i.e. is a Lie symmetry of $({\cal S}_{\Gamma})$. So the study of $Aut(\Gamma)$ can be reduced as a very special case to the general problem of 
study of Lie symmetries of holomorphic completely overdetermined second 
order involutive systems.  We emphasize that {\it the systems describing the Segre families 
form a very special subclass in the class of holomorphic completely overdetermined second 
order involutive systems with one dependent variable}. So  theorem \ref{theo3.2} generalize 
known results in several directions.

Indeed, since $Sym({\cal S}_{\Gamma})$ is a complex Lie transformation group in view of 
 our theorem and $Aut(\Gamma)$ is its closed subgroup, we conclude that 
it is a local real Lie transformation subgroup of $Sym({\cal S}_{\Gamma})$. 
Hence   theorem \ref{theo3.2} (which we apply in the special case 
of one dependent variable) gives an upper estimate for its dimension: the real 
dimension of $Aut(\Gamma)$ is majorated by $2dim Sym({\cal S}_{\Gamma})$, in 
particular, by $2(n^2 + 4n + 3)$. In order to improve this estimate, we recall the following useful observation due to 
E.Cartan \cite{Ca}. Let a holomorphic vector field $X$  generate 
a local real one-dimensional subgroup of $Aut(\Gamma)$. This is 
equivalent to the fact that $Re X$ is a tangent vector field to $\Gamma$. 
 Since  this subgroup is a real one-parameter subgroup of $Sym({\cal S}_{\Gamma})$, 
we have necessarily  $X \in Lie({\cal S}_{\Gamma})$. 
  $\Gamma$ is Levi nondegenerate, so the field $Re(iX)$ cannot be 
tangent to $\Gamma$ simultaneously with $Re X$ i.e.  $Lie(\Gamma)$
is a totally real subspace of $Lie({\cal S}_{\Gamma})$. Therefore, the real 
dimension of $Aut(\Gamma)$ is majorated by  the complex dimension 
of $Lie({\cal S}_{\Gamma})$. In particular, it is smaller that $n^2 + 4n + 3$. The example of the sphere shows 
that this estimate is precise.  We obtain the following

\begin{e-cor}
  $Aut(\Gamma)$ is a local real Lie transformation subgroup of 
$Sym({\cal S}_{\Gamma})$ and the dimension of $Aut(\Gamma)$ is majorated by 
the complex dimension of $Sym({\cal S}_{\Gamma})$. In particular, it is always 
majorated by $n^2 + 4n + 3$. Moreover,  every infinitesimal automorphism of $\Gamma$ is uniquely determined 
by its second order Taylor developement at a fixed point.
\end{e-cor}

Thus, we find here some classical results of N.Tanaka \cite{Ta},
S.S.Chern - J.Moser \cite{CM} in the infinitesimal form.

\bigskip

{\bf 2. Lie method and proof of the main theorem}

\bigskip

 Let $G$ be a local group of biholomorphic transformation acting on a domain in 
$\cc^n \times \cc^m$. Every  biholomorphism $g \in G$ ,  
$g: (x,u) \mapsto (x^*,u^*)$ {\it close enough to the identity}   
lifts canonically to a fiber preserving biholomorphism $g^{(r)}: J^r_{n,m} 
\longrightarrow J^r_{n,m}$ as follows: if $u = f(x)$ is a holomorphic function 
near $p$, $q = f(p)$ and $u^* = f^*(x^*)$ is its image under $g$ 
(i.e. the graph of $f^*$ is the image of the graph of $f$ under $g$ near 
the point $(p^*, q^*) = g(p,q)$), then the jet $j^r_{p^*}(f^*)$ is by the definition the image of $j^r_p(f)$ under 
$g^{(r)}$. In particular, a one-parameter local Lie group of 
transformations $G$ canonically lifts to $J^r_{n,m}$ as a one-parameter Lie group 
of transformations $G^{(r)}$  which is called the {\it $r$-prolongation} of $G$. The infinitesimal
generator $X^{(r)}$ of $G$ is called the {\it $r$-prolongation} of the infinitesimal 
generator $X$ of $G$. 

Let $({\cal S})$ be a holomorphic PDE system of $r$th order with $n$ independent and $m$ 
dependent variables. Then it defines naturally a complex subvariety $({\cal S})^{(r)}$ in the jet 
space $J^{(r)}_{n,m}$ obtained by  replacing of 
 derivatives of  dependent variables  by the corresponding coordinates in the
jet space.  So $u$ is a holomorphic solution of the system $({\cal S})$ if and only 
if the section $(p,u(p),j^r_p(u))$ of the holomorphic fibre bundle 
$\pi_n: J^r_{n,m} \longrightarrow \cc^n$ (with the natural projection $\pi_n$) 
is contained in the 
variety $({\cal S})^{(r)}$. A key proposition of the Lie theory states  that if 
the r-prolongation $X^{(r)}$ of a vector field $X$ 
is a tangent field to  $({\cal S})^{(r)}$ then $X$ 
is an infinitesimal symmetry of $({\cal S})$  (\cite{Ol}, Theorem 2.31). 
 
For the system (\ref{1}) one has  

\begin{eqnarray*}
({\cal S})^{(2)}: u^k_{ij} = \hat F_{ij}^k(x,u,u^{(1)})
\end{eqnarray*}
Here $\hat F^{k}_{ij}$ denote the natural lifting of $F^{k}_{ij}$ to 
the jet space obtained via the replacing of the derivatives $u^{\mu}_{x_i}$ by 
the jet coordinates $u^{\mu}_i$ in power expansions of $F^{k}_{ij}$, i.e. 
$\hat F^{k}_{ij} = F^{k}_{ij}(x,u,u^{(1)})$.

In the natural coordinates one has 
\begin{eqnarray*}
X^{(r)} = X + \sum_{j,\mu} \eta^{\mu}_j\frac{\partial}{\partial u^{\mu}_j} +...
+ \sum_{i_1,...,i_r,\mu} \eta^{\mu}_{i_1i_2...i_r}\frac
{\partial}{\partial u^{\mu}_{i_1i_2...i_r}}
\end{eqnarray*}
 where $\eta^{\mu}_i = D_i\eta^{\mu} - \sum_j (D_i\theta_j) u^{\mu}_j$, 
$\eta^{\mu}_{i_1...i_{r-1}i_r} = D_{i_r} \eta^{\mu}_{i_1...i_{r-1}} - 
\sum_j (D_{i_r}\theta_j)u^{\mu}_{i_1...i_{r-1}j}$ and 
\begin{eqnarray*}
D_i = \frac{\partial}{\partial x_i} + \sum_k u^k_i \frac{\partial}{\partial u^k} + 
\sum_{\mu,j} u^{\mu}_{ij}\frac{\partial}{\partial u^{\mu}_j} + ... 
\end{eqnarray*}
is the operator of {\it total derivative} \cite{BlKu}, \cite{Ol}.

In our case  a direct computation  gives an explicit expression 
for the coefficients of $X^{(2)}$:

\begin{eqnarray*}
& &\eta^{\mu}_{i_1} = \frac{\partial \eta^{\mu}}{\partial x_{i_1}} + \sum_k u^k_{i_1}
\frac{\partial \eta^{\mu}}{\partial u^k} - \sum_j \left ( \frac{\partial \theta_j}
{\partial x_{i_1}} + \sum_k u^k_{i_1} \frac{\partial \theta_j}{\partial u^k} \right ) 
u^{\mu}_j ,\\
& &\eta^{\mu}_{i_1i_2} = \frac{\partial^2\eta^{\mu}}{\partial x_{i_2} 
\partial x_{i_1}} +u^{\mu}_{i_1} \left [\frac{\partial^2\eta^{\mu}}
{\partial x_{i_2}\partial u^{\mu}} - \frac{\partial^{2}\theta_{i_1}}
{\partial x_{i_2}\partial x_{i_1}} \right ] + u^{\mu}_{i_2}
\left [\frac{\partial^2\eta^{\mu}}{\partial x_{i_1}\partial u^{\mu}}Ê- 
\frac{\partial^2\theta_{i_2}}{\partial x_{i_2}\partial x_{i_1}}\right ]\\
& &+ \sum_{k \neq \mu} u^k_{i_1}\frac{\partial^2\eta^{\mu}}
{\partial x_{i_2}\partial u^k}+ \sum_{k\neq \mu} u^k_{i_2}\frac{\partial^2\eta^{\mu}}
{\partial x_{i_1}\partial u^k} - 
\sum_{k\neq i_1,k\neq i_2}u^{\mu}_k\frac{\partial^2\theta_k}
{\partial x_{i_2}\partial x_{i_1}}- 
\sum_{k; j\neq i_2}u^k_{i_1}u^{\mu}_{j}\frac{\partial^2\theta_j}
{\partial x_{i_2}\partial u^k}\\
& &- \sum_{i; s\neq i_1} u^i_{i_2}u^{\mu}_{s}\frac{\partial^2\theta_s}
{\partial x_{i_1}\partial u^i} + \sum_{r\neq \mu, p\neq \mu}u^r_{i_2}u^{p}_{i_1}
\frac{\partial^2 \eta^{\mu}}{\partial u^r \partial u^p} + 
\sum_{t \neq \mu} u^t_{i_1}u^{\mu}_{i_2}\left [-\frac{\partial^2\theta_{i_2}}{\partial x_{i_2}\partial u^t} + 
\frac{\partial^2\eta^{\mu}}{\partial u^{\mu} \partial u^t}\right ]\\
& & + \sum_{q \neq \mu} u^q_{i_2}u^{\mu}_{i_1}\left [-\frac{\partial^2\theta_{i_1}}
{\partial u^q \partial x_{i_1}} + \frac{\partial^2 \eta^{\mu}}
{\partial u^q \partial u^{\mu}}\right ] + \left [\frac{\partial^2\eta^{\mu}}
{(\partial u^{\mu})^2} - \frac{\partial^2\theta_{i_2}}
{\partial x_{i_2}\partial u^{\mu}} - \frac{\partial^2 \theta_{i_1}}
{\partial x_{i_1} \partial u^{\mu}} \right ]
u^{\mu}_{i_1}u^{\mu}_{i_2}\\
& & - 
\sum_{a,b,s} u^a_{i_2}u^b_{i_1}u^{\mu}_{s}\frac{\partial^2\theta_s}
{\partial u^a \partial u^b} + \Lambda^{\mu}_{i_1i_2}
\end{eqnarray*}
for $i_1 \neq i_2$ and

\begin{eqnarray*}
& & \eta^{\mu}_{ii} = \frac{\partial^2\eta^{\mu}}{(\partial x_i)^2} + 
u^{\mu}_i\left [2\frac{\partial^2\eta^{\mu}}
{\partial x_i\partial u^{\mu}} - \frac{\partial^2 \theta_{i}}{(\partial x_i)^2}\right ] + 
2\sum_{k \neq \mu}u^k_i\frac{\partial^2\eta^{\mu}}{\partial x_i\partial u^k}  - 
\sum_{k \neq i}u^{\mu}_k\frac{\partial^2\theta_k}{(\partial x_i)^2}\\
& & - 
2\sum_{k; j \neq i}
u^k_i u^{\mu}_j\frac{\partial^2\theta_j}{\partial x_i \partial u^k}+ 
\sum_{r\neq \mu; p\neq \mu}
u^r_i u^p_i\frac{\partial^2\eta^{\mu}}{\partial u^r\partial u^p} + 
\sum_{t \neq \mu}u^t_iu^{\mu}_i\left [- \frac{\partial^2\theta_i}
{\partial x_i \partial u^t} +
\frac{\partial^2\eta^{\mu}}{\partial u^{\mu} \partial u^t}\right ]+\\
& & \sum_{q \neq \mu} u^q_iu^{\mu}_i
\left [-\frac{\partial^2\theta_i}{\partial x_i \partial u^q} + 
\frac{\partial^2\eta^{\mu}}
{\partial u^q \partial u^{\mu}}\right ] + 
\left [\frac{\partial^2\eta^{\mu}}{(\partial u^{\mu})^2} - 
2 \frac{\partial^2\theta_i}
{\partial x_i \partial u^{\mu}}\right ](u^{\mu}_i)^2 \\
& &- \sum_{a,b,s}u^a_iu^b_iu^{\mu}_s\frac
{\partial^2 \theta_s}{\partial u^a \partial u^b} + \Lambda^{\mu}_{ii}
\end{eqnarray*}
with

\begin{eqnarray*}
& &\Lambda^{\mu}_{i_1i_2} = \sum_s u^s_{i_2i_1} \frac{\partial \eta^{\mu}}{\partial u^s} -
\sum_p u^{\mu}_{i_2p} \frac{\partial \theta_p}{\partial x_{i_1}} - \sum_j u^{\mu}_{i_1j}
\frac{\partial \theta_j}{\partial x_{i_2}} - 
\sum_{p,q} u^q_{i_2i_1} u^{\mu}_p \frac{\partial \theta_p}{\partial u^q}\\
& & - \sum_{p,q} u^{\mu}_{i_2p}u^q_{i_1} \frac{\partial \theta_p}{\partial u^q} - 
\sum_{s,j} u^{\mu}_{i_1j}u^s_{i_2}\frac{\partial \theta_j}{\partial u^s}
 \end{eqnarray*}

 Since the system (\ref{1}) is involutive,  
for every point $P \in ({\cal S})^{(2)}$ with 
the natural projection $\pi_{n,m}(P) = (p,q) \in \cc^n \times \cc^m$ 
there exists a solution $u(x)$ of $({\cal S})$ holomorphic near $p$ such that 
$(p,q,j^2_p(u)) = P$. 
So the Lie criterion  (\cite{Ol}, Theorem 2.72) implies  that 
$X \in Lie({\cal S})$ if and only if the second prolongation 
satisfies the following system of equations in $J^2_{n,m}$: 

\begin{eqnarray*}
X^{(2)}(u^{\mu}_{ij} - 
\hat F^{\mu}_{ij}) = 0, u^{\mu}_{ij} = \hat F^{\mu}_{ij}
\end{eqnarray*}

 This system implies that 
$\eta^{\mu}_{ij} = X^{(2)}(\hat F^{\mu}_{ij}) = X^{(1)}(\hat F^{\mu}_{ij})$. 
Replace now in the expressions of 
$\Lambda^{\mu}_{i_1i_2}$ the jet coordinates $u^{k}_{ij}$ by $\hat F^{k}_{ij}$ and 
denote obtained functions by $\hat \Lambda^{\mu}_{ij}$. Transfer them to the right 
side; we get the equations of the form 

\begin{eqnarray}
\label{5}
\tilde\eta^{\mu}_{ij} = X^{(1)}(\hat F^{\mu}_{ij}) - \hat \Lambda^{\mu}_{ij}
\end{eqnarray}  

Without loss of generality assume that every $\hat F^{\mu}_{ij}$ is
  represented by a power series with respect to $u^{k}_s$. 
Then we can  develop the right sides 
$X^{(1)}(\hat F^{\mu}_{ij}) - \hat \Lambda^{\mu}_{ij}$ of our equations 
in power series with respect to $u^k_l$.  Clearly,  the coefficients of these 
expansions 
are completely determined by the coefficients 
of the expansions of $\hat F^{\mu}_{ij}$  and one can  effectively compute them 
in a concrete case.

Comparing now the coefficients near the powers of $u^{k}_l$ of degree $\leq 3$ in 
the equations (\ref{5}) and using 
explicit expressions for the coefficients of $X^{(2)}$, we obtain the 
following PDE system for the coefficients of the inifinitesimal symmetry $X$:

\begin{eqnarray*}
& &(A): \frac{\partial^2\eta^{\mu}}{\partial x_{i_2}\partial x_{i_1}} = A^{\mu}_{i_1i_2},
\frac{\partial^2\eta^{\mu}}{\partial x_i \partial u^k} = B^{\mu}_{ik}, k \neq \mu, 
\frac{\partial^2\eta^{\mu}}{\partial u^r \partial u^p} = C^{\mu}_{rp}, r\neq \mu, p\neq \mu,\\
& &(B): \frac{\partial^2\theta_k}{\partial x_{i_2}\partial x_{i_1}} = D^k_{i_1i_2}, 
k \neq i_1, k\neq i_2, \frac{\partial \theta_j}{\partial x_i \partial u^k} = E^j_{ik}, j \neq i, (C): \frac{\partial^2 \theta_s}{\partial u^a \partial u^b} = G^{s}_{ab}\\
& & (D_1): \frac{\partial^2\eta^{\mu}}{\partial x_{i_2}\partial u^{\mu}} - 
\frac{\partial^2 \theta_{i_1}}{\partial x_{i_2}\partial x_{i_1}} = 
H^{\mu}_{i_1i_2}, i_1 \neq i_2, (D_2): \frac{\partial^2\eta^{\mu}}{\partial u^{\mu} \partial u^t} - \frac{\partial^2 
\theta_{i_2}}{\partial x_{i_2} \partial u^t} = I^{\mu}_{i_2t}, t\neq \mu\\
& & (D_3): 2\frac{\partial^2\eta^{\mu}}{\partial x_i \partial u^{\mu}} - \frac
{\partial^2 \theta_i}{(\partial x_i)^2} = J^{\mu}_i, (D_4): \frac{\partial^2\eta^{\mu}}{(\partial u^{\mu})^2} - 2\frac{\partial^2 \theta_i}
{\partial x_i \partial u^{\mu}} = K^{\mu}_i
\end{eqnarray*}
where  the right sides are analytic functions in 
 $\lambda(x,u) = (x,u, \alpha(x,u),\beta(x,u),\delta(x,u),\varepsilon(x,u))$.  
 We point out that we do not write here all obtained equations; we consider only 
those who will be enough for the proof of our results.

Denote by $\Omega$ a holomorphic vectorvalued function whose components coincide with the right 
sides of our system: $\Omega = (A^{\mu}_{i_1i_2},B^{\mu}_{ik},..., K^{\mu}_i)$.

 An important propery of obtained PDE system $(A)- (D_i)$ is 
its {\it linearity} with respect to the second order derivatives of dependent 
variables.  Denote by $v$ the vector $\cc^L$ (for a suitable $L$) 
whose components are the second order partial derivatives of $\theta_j $ and 
$\eta^{\mu}$; then our system can be written in the form $M v = \Omega$,
where $M$ is an integer matrix.  An elementary linear algebra argument shows that 
this system can be represented  in  the form $M' v' = P \gamma + \Omega$, 
where $v'$ is a vector formed by components of $v$ which are not components of 
$\gamma$, $P$ is an integer  matrix  and $M'$ is an {\it invertible} 
square integer matrix.  Therefore 
$v' = M^{-1}P \gamma + M^{-1}\Omega$, so  every second order partial derivative 
of $\theta_j$, $\eta^k$ is a linear combination of components of $\gamma$, $\Omega$.
In particular, {\it the second order partial derivatives of $\theta_j$, $\eta^k$ 
at $(x,u)$ are determined by $\omega(x,u)$}.   
 Denote by $V$ the vector whose components are the third order partial 
derivatives of $\theta_j$, $\eta^{\mu}$ and write the system obtained by taking the 
partial derivatives in $(A) - (D_i')$ in the form $N V =  \Omega'$
where $ \Omega'$ denote the vector with components 
$\frac{\partial\Omega_j(\lambda(x,u))}{\partial x_i}$, 
$\frac{\partial \Omega_j(\lambda(x,u))}{\partial u^k}$ and $N$ is an  integer matrix. 
A direct computation shows  that one can choose a subsystem of this 
system with an invertible square matrix $N'$, so we obtain that 
for every multi-indice $\tau$, $\vert \tau \vert = 3$, there are  polynomials
 $R^{\tau}_k$, $S^{\tau}_j$  with rational coefficients, 
such that the following holds: 
\begin{eqnarray}
\label{6}
\frac{\partial^{3} \theta_j}{\partial x_1^{\tau_1}...\partial x_n^{\tau_n}
\partial (u^1)^{\tau_{n+1}}...\partial (u^m)^{\tau_{n+m}}}(x,u) = 
R^{\tau}_j((\partial \Omega)
(\lambda(x,u)),(\partial^{2}\theta)(x,u),(\partial^{2}\eta)(x,u))
\end{eqnarray}
\begin{eqnarray}
\label{7}
 \frac{\partial^{3} \eta^k}{\partial x_1^{\tau_1}...\partial x_n^{\tau_n}
\partial (u^1)^{\tau_{n+1}}...\partial (u^m)^{\tau_{n+m}}}(x,u) = 
S^{\tau}_k((\partial \Omega)(\lambda(x,u)),
(\partial^{2}\theta)(x,u),(\partial^{2}\eta)(x,u))
\end{eqnarray}
where $(\partial \Omega)(\lambda(x,u))$ denote the vector function 
 whose components are the first order partial derivatives of $\Omega_j$ evaluated at $\lambda(x,u)$, 
$(\partial^2 \theta)(x,u)$ (resp. $(\partial^2 \eta)(x,u)$) denote the vector 
function whose components are the partial derivatives of all 
$\theta_j$ (resp. $\eta^k$) of order $\leq 2$.  
This means that {\it the third order partial derivatives of $\theta_j$, $\eta^k$ 
at $(x,u)$ are determined by $\omega(x,u)$}. 
If $\omega(x,u)$ is given now, (\ref{6}), (\ref{7}) and the chain rule show
 that all coefficients of the Taylor expansions of $\theta_j$, $\eta^k$ are 
determined by recursion. This completes the proof of the theorem.

It is also clear from the construction that the corresponding 
homogeneous system $M v = 0$ describes infinitesimal symmetries of the  
system $({\cal S}_0)$ of the form (\ref{1}) with $F^k_{ij} \equiv 0$.  We obtain that 
the set $Lie({\cal S}_{0})$ of infinitesimal symmetries of this 
system  is a complex Lie algebra of dimension $(n + m + 2)(n + m)$ 
generated by the following holomorphic vector fields: 
$U_k = \frac{\partial}{\partial x_k}$, $V_{\mu} = \frac{\partial}{\partial u^{\mu}}$, 
$W_{jk} = x_j\frac{\partial}{\partial x_k}$, $A_{jk} = u^j\frac{\partial}{\partial x_k}$
, $B_{k\mu} = x_k\frac{\partial}
{\partial u^{\mu}}$, $C_{k\mu} = u^k\frac{\partial}{\partial u^{\mu}}$,
$X_j = \sum_k x_jx_k \frac{\partial}{\partial x_k} + \sum_{mu}x_ju^{\mu}\frac
{\partial}{\partial u^{\mu}}$, $Y_{\nu} = \sum_k x_ku^{\nu} \frac{\partial}{\partial x_k}
+ \sum_{mu} u^{\nu}u^{\mu}\frac{\partial}{\partial u^{\mu}}$, $k,j = 1,...,n$, 
$\mu,\nu = 1,...,m$. In particular, $dim Lie({\cal S}_0) = (n + m + 2)(n+m)$ 
so the dimension estimate given by our theorem is precise.

In the present paper we restrict an application of the Lie method 
only by the classical case of a Levi nondegenerate hypersurface. But this method 
can also be applied in other situations which are of interest for the CR 
geometry and  form an area of research activity of several authors.
 For instance, the Levi degenerate case leads to considerations of holomorphic 
second order completely overdetermined involutive PDE systems which are not 
solved with respect to the second order derivatives. A study of their symmetries 
requires a combination of the Lie method with  some tools of the local complex 
analytic geometry (compare with \cite{DW}, \cite{DF}). On the other hand, the Segre families of Cauchy - Riemann 
manifolds of higher codimension are described by  holomorphic second order 
completely overdetermined involutive PDE systems with additional first order 
relations, i.e. with additional holomorphic equations including first order derivatives of 
dependent variables.   Clearly, the Lie method allows to study this class of systems and
  just requires more involved computations.

\end{document}